V.I. NORKIN, A.Y. KOZYRIEV

# ON SHOR'S R-ALGORITHM FOR PROBLEMS WITH CONSTRAINTS

**Abstract.** *Shor's r-algorithm (Shor, Zhurbenko (1971), Shor (1979)) with space stretching in the direction of difference of two adjacent subgradients is a competitive method of nonsmooth optimization. However, the original r-algorithm is designed to minimize convex ravine functions without constraints. The standard technique for solving constraint problems with this algorithm is to use exact nonsmooth penalty functions (Eremin (1967), Zangwill (1967)). At the same time, it is necessary to choose the (sufficiently large) penalty coefficient in this method. In Norkin (2020, 2022) and Galvan et al. (2021), the so-called projective exact penalty function method is proposed, which does not formally require an exact definition of the penalty coefficient. In this paper, a nonsmooth optimization problem with convex constraints is first transformed into a constraint-free problem by the projective penalty function method, and then the r-algorithm is applied to solve the transformed problem. We present the results of testing this approach on problems with linear constraints using a program implemented in Matlab.*

**Introduction**. Nonsmooth optimization problems arise in a wide range of relevant fields, including engineering, finance and deep learning, where mostly chosen activation functions have discontinuous derivatives, e.g., ReLU. Conventional optimization algorithms, designed predominantly for smooth problems, face difficulties applied to nonsmooth contexts due to discontinuities and other associated irregularities.

One approach to overcome these difficulties is to apply smoothing of functions [1 – 3]. The other alternative is to apply methods of nondifferentiable optimization [4, 5]. For both approaches the question is how to treat constraints. In the paper we investigate exact penalty methods [6 – 11] for this purpose. As a nonsmooth optimizer we explore Shor's r-algorithm [4, 12, 13], which is a competitive non-smooth optimization method [5]. However, the original method is designed to minimize convex functions without constraints. A standard technique for solving constrained problems by this algorithm is to apply exact nonsmooth penalty functions [6 – 8]. At the same time, it is necessary to choose correctly the (sufficiently large) penalty coefficient in the penalty function method. In the works of [9 – 11] authors proposed the so-called exact projective penalty method, which theoretically does not require the selection of the penalty coefficient.

In this paper, a nonsmooth optimization problem with convex constraints is first transformed into an unconstrained problem by the method of projective penalty functions, and then the r-algorithm is applied to solve the transformed problem. We provide testing results of this approach for an optimization problem with linear constraints using Matlab program for the r-algorithm.

**The problem setting.** Suppose that the problem of convex optimization is solved:

$$f(x) \to \min_{x\in X}\,, \qquad (1)$$

where $f(x)$ is a convex (nonsmooth) objective function on a convex set $X \subset \mathbb{R}^n$; $\mathbb{R}^n$ is an $n$-dimensional arithmetic Euclidean space.

The classical penalty function method for the problem with constraints $x \in X = \{x \in \mathbb{R}^n : h(x) \le 0\}$ and convex function $h(x)$ has the form [6 – 8]:

$$F(x) = f(x) + M \max\{0, h(x)\} \to \min_{x \in \mathbb{R}^m}, \qquad (2)$$

where $M > 0$ is a sufficiently large preliminary unknown constant. Another variant of the classical exact penalty method has the form [16, Proposition 9.68]

$$F(x) = f(x) + M \|x - \pi_X(x)\|^\gamma \to \min_{x \in \mathbb{R}^m}, \qquad 0 < \gamma \le 1, \qquad (3)$$

where $\pi_X(x) = \arg\min_{y \in X} \|y - x\|^2$ is the Euclidean projection of an element $x \in R^n$ into a convex set $X$; for $x \in X$, obviously, $\pi_X(x) = x$. Here again $M$ must be sufficiently large for Lipschitzian $f$ and $\gamma = 1$, and $F$ is non-Lipshitzian when $\gamma \in (0,1)$. If the set $X$ is defined by linear constraints, then finding $\pi_X(x)$ is a quadratic programming problem. With simple restrictions, for example, with two-sided restrictions on variables or restrictions in the form of a simplex, this problem can be solved analytically. Other valid (non-Euclidean) projections are considered in [9, 10, 15].

The authors of [9 – 11] proved that a complex continuous problem (1) with convex constraint set $X$ can be equivalently reduced to an unconditional optimization problem with penalty function of the form:

$$F(x) = f(\pi_X(x)) + M \|x - \pi_X(x)\|^\gamma \to \min_{x \in \mathbb{R}^m}, \qquad (4)$$

where $\gamma > 0$ and $M > 0$ are arbitrary positive numbers.

The function $F(x)$ in (4) is nonsmooth, and at $\gamma = 1$ it is Lipschitzian. It may be non-convex, but it does not have local minima different from the solution of problem (1) [10, 11]. The advantage of form (4) with respect to (2) and (3) is that penalty parameters $M$ and $\gamma$ in (4) can be arbitrary chosen but in (2) $M$ has to be carefully selected.

**r-Algorithm.** In this paper we reduce the original constrained problem (1) to unconstrained ones (3) and (4) and solve the latter by Shor's r-algorithm [4, 12]. The principal r-algorithm and its implementations are described in detail in a number of works [4, 12 – 14, 17], so we do not provide here its conceptual description. We just adapted the Octave version of the r-algorithm from [13, 14] to Matlab system [18]. For solving problem (3) with $\gamma = 1$, we use the following subgradient $g_F(x)$ of the penalty function $F(x)$, $g_F(x) = g_f(x) + M\left(x - \pi_X(x)\right)/\|x - \pi_X(x)\|$, where $g_f(x)$ is a subgradient of $f(\cdot)$ at $x$ and 0/0 is considered as $0$. For solving problem (4) we do not use a program for calculating subgradients of the function but heuristically use finite-difference estimates of subgradients of the penalty function $F(x)$ (the finite difference step is equal to $\sqrt{eps} = 1.4901\text{e-}08$, which is the default value for finite differences in Matlab, $eps$ is the smallest positive number in the Matlab system [18]). The penalty factor $M$ in (4) can theoretically be arbitrary; in numerical experiments it varied within the limits $[10^1, 10^4]$, the parameter $\gamma$ was 1. Numerical experiments were carried out on the separable ravine function

$$f(x) = \sum_{i=1}^{n} (1.2)^{i-1} |x_i - 1|, \qquad n \in [10, 100], \qquad (5)$$

under box constraints $X=\left\{x\in\mathbb{R}^n : x_i\in\left[c_i,d_i\right], i=1,...,n\right\}$. The solution of the problem (5) is: $x_i^*\in\max\left\{c_i,\min\left\{1,d_i\right\}\right\}$, $i=1,...,n$. If $1\in\left[c_i,d_i\right]$, then $x^*=\left(1,1,...,1\right)^T$, $f^*=f^*(x^*)=0$. In this case the projection $\pi_X(x)$ of point $x$ on $X$ is found analytically, $\left(\pi_X(x)\right)_i=\max\left\{c_i,\min\left\{x_i,d_i\right\}\right\}$, $i=1,...,n$.

Besides, we consider problem (5) under linear constraints

$$X=\left\{Ax\le b,\quad A_{eq}x=b_{eq},\quad b\in\mathbb{R}^{m_1},\quad b_{eq}\in\mathbb{R}^{m_2},\quad x_i\in\left[c_i,d_i\right],\quad i=1,...,n.\right\}.$$

In this case, the projection $\pi_X(x)$ is found by means of Matlab specialized quadratic programming solver ***lsqlin***. In numerical experiments we used a particular feasible set:

$$X=\left\{\sum\nolimits_{i=1}^{n}x_i\le b,\quad b\in[0,n]\quad x_i\in\left[0,1\right],\quad i=1,...,n\right\}. \tag{6}$$

The optimal solution $x^*$ of problem (5)-(6) is the following: $x_i^*=0$, $i=1,...,n-\lfloor b\rfloor-1$; $x_{n-\lfloor b\rfloor}^*=b-\lfloor b\rfloor$; $x_i^*=1$, $i=n-\lfloor b\rfloor+1,...,n$; $\lfloor b\rfloor$ is the integer part of $b$, i.e. the largest integer smaller or equal to $b$. The accuracy of an approximate solution $x$ is measured as $\delta=\max_{1\le i\le n}\left|x_i-x_i^*\right|$, $\varepsilon=\left|f(x)-f(x^*)\right|$.

The parameters of the r-algorithm are as follows.

```
%r_algorithm_data, version ralgb5
h0=norm(ub-lb);%Initial step
alpha=4; q1=1; q2=1.1; nh=3;%Shor's algorithm parameters
epsx=10^-8; epsg=10^-12; maxitn=7000;%Shor's algorithm parameters
```

**Numerical experiments**. Numerical experiments were conducted on HP personal computer 11th Gen Intel(R) Core(TM) i5-1135G7 @ 2.40GHz 16Gb RAM. The results of numerical experiments are presented in Tables 1 – 4, where parameter $n$ designates the dimension of the space, $M$ is the penalty coefficient. Cells contain information on the accuracy $\varepsilon$ or $\delta$ , number of r-algorithm iterations *itn*, and solving time in seconds.

Table 1 contains results of solving problem (5) by means of penalty method (3) under box constraints: $x_i\in\left[0,1\right]$, $i=1,...,n$. In this example the value of $M$ does not play a significant role in calculations, and computational times are small because the problem is separable.

TABLE 1. Testing the exact penalty method (3) (under box constraints: $x_i\in\left[0,1\right]$, $i=1,...,n$ ).

| | n = 10 | n = 20 | n = 30 | n = 50 | n = 100 |
|---|---|---|---|---|---|
| M = 1 | $\varepsilon$ = 0.000000e+00<br>itn = 117<br>time sec = 0.0732 | $\varepsilon$ = 0.000000e+00<br>itn = 242<br>time sec = 0.0956 | $\varepsilon$ = 0.000000e+00<br>itn= 381<br>time sec = 0.1336 | $\varepsilon$ = 0.000000e+00<br>itn = 720<br>time sec = 0.3140 | $\varepsilon$ = 0.000000e+00<br>itn = 1704<br>time sec = 0.6013 |
| M = 10000 | $\varepsilon$ = 0.000000e+00<br>itn = 199<br>time sec = 0.1904 | $\varepsilon$ = 0.000000e+00<br>itn = 1110<br>time sec = 0.9118 | $\varepsilon$ = 0.000000e+00<br>itn = 3579<br>time sec = 2.5173 | $\varepsilon$ = 0.000000e+00<br>itn = 1105<br>time sec = 0.5193 | $\varepsilon$ = 0.000000e+00<br>itn = 2420<br>time sec = 1.0253 |

Table 2 contains results of solving problem (5) by means of penalty method (3) under constraints: $\sum\nolimits_{i=1}^{n}x_i\le n/2$, $x_i\in\left[0,1\right]$, $i=1,...,n$. In this case the value of $M$ is important and has to be chosen large enough. Empty cells mean that r-algorithm stopped at an infeasible point.

TABLE 2. Testing the exact penalty method (3) (under constraints: $\sum_{i=1}^{n} x_i \le n/2$, $x_i \in [0,1]$, $i = 1,...,n$ ).

| | n =10 | n = 20 | n = 30 | n = 40 | n = 50 |
|---|---|---|---|---|---|
| M = 1 | – | – | – | – | – |
| M = 10 | δ = 7.631239e-04<br>itn = 1000<br>time sec = 5.6612 | – | – | – | – |
| M = 100 | δ = 6.821032e-04<br>itn = 1000<br>time sec = 4.7020 | δ = 3.375257e-03<br>itn = 2000<br>time sec = 6.5575 | δ = 3.508656e-04<br>Itn=487<br>time sec = 0.7939 | – | – |
| M = 1000 | δ = 6.697289e-04<br>itn = 2000<br>time sec = 9.6003 | δ = 8.498720e-04<br>itn = 2000<br>time sec = 6.3643 | δ = 1.854921e-03<br>itn = 2000<br>time sec = 6.3189 | δ = 2.265069e-03<br>itn = 673<br>time sec = 2.8188 | – |
| M = 10000 | δ = 6.454530e-04<br>itn = 1000<br>time sec = 4.7005 | δ = 6.534425e-04<br>itn = 1000<br>time sec = 3.1087 | δ = 2.082616e-03<br>itn = 1000<br>time sec = 2.6072 | δ = 3.839708e-03<br>itn = 1000<br>time sec = 2.4281 | δ = 7.590022e-03<br>itn = 804<br>time sec= 1.2983 |

Table 3, 4 contain results of solving problem (5) by means of penalty method (4) under box constraints: $x_i \in [0,1]$, $i = 1,...,n$, and under constraints $\sum_{i=1}^{n} x_i \le n/2$, $x_i \in [0,1]$, $i = 1,...,n$, respectively. In these examples the value of $M$ does not play any significant role in calculations but the execution times are large compared to the ones of Table 1, 2. This is explained by time consuming finite-difference estimation of generalized gradients of the penalty function (4). If finite differences were calculated parallelly, then execution times would potentially be reduced by $n$ and be comparable with the ones from Tables 1, 2.

TABLE 3. Testing the exact projective penalty function method (4) (box constraints: $x_i \in [0,1]$, $i = 1,...,n$ ).

| | n = 10 | n = 20 | n = 30 | n = 50 | n = 80 |
|---|---|---|---|---|---|
| M = 1 | $\varepsilon$ = 3.296448e-08<br>itn = 163<br>time sec = 2.6436 | $\varepsilon$ = 7.550370e-08<br>itn = 352<br>time sec = 6.1961 | $\varepsilon$ = 1.049490e-07<br>itn = 605<br>time sec = 15.6210 | $\varepsilon$ = 2.720405e-07<br>itn = 1647<br>time sec = 74.2272 | – |
| M = 10 | $\varepsilon$ = 4.257959e-09<br>itn = 148<br>time sec = 1.0699 | $\varepsilon$ = 2.944758e-07<br>itn = 293<br>time sec = 3.8003 | $\varepsilon$ = 5.839912e-07<br>itn = 445<br>time sec = 8.8346 | $\varepsilon$ = 1.381021e-06<br>itn = 902<br>time sec = 33.0274 | $\varepsilon$ = 3.652280e-06<br>itn = 2417<br>time sec = 245.577 |
| M = 100 | $\varepsilon$ = 1.516709e-08<br>itn = 214<br>time sec = 1.4820 | $\varepsilon$ = 1.579981e-08<br>itn = 322<br>time sec = 4.2364 | $\varepsilon$ = 2.446992e-06<br>itn = 441<br>time sec = 7.9904 | $\varepsilon$ = 7.255148e-06<br>itn = 807<br>time sec = 26.6567 | $\varepsilon$ = 2.555098e-05<br>itn = 1667<br>time sec = 148.264 |
| M = 1000 | $\varepsilon$ = 1.330437e-10<br>itn = 225<br>time sec = 1.6897 | $\varepsilon$ = 2.487261e-08<br>itn = 595<br>time sec = 9.2577 | $\varepsilon$ = 1.143795e-06<br>itn = 549<br>time sec = 10.2455 | $\varepsilon$ = 4.691534e-05<br>itn = 791<br>time sec = 24.1373 | $\varepsilon$ = 9.800041e-05<br>itn = 1572<br>time sec = 129.499 |
| M = 10000 | $\varepsilon$ = 5.581399e-09<br>itn = 215<br>time sec = 1.7683 | $\varepsilon$ = 6.430452e-07<br>itn = 1422<br>time sec = 35.8907 | $\varepsilon$ = 4.211114e-03<br>itn = 5000<br>time sec = 175.070 | $\varepsilon$ = 2.912054e-06<br>itn = 1060<br>time sec = 36.3398 | $\varepsilon$ = 8.120949e-04<br>itn = 1427<br>time sec = 125.728 |

TABLE 4. Testing the exact penalty method (4) (under constraints: $\sum_{i=1}^{n} x_i \le n/2$, $x_i \in [0,1]$, $i = 1,...,n$ ).

| | n = 10 | n = 20 | n = 30 | n = 40 | n=80 |
|---|---|---|---|---|---|
| M = 1 | $\delta = 2.888739e-04$<br>itn = 157<br>time sec = 4.4791 | $\delta = 5.075395e-04$<br>itn = 339<br>time sec = 17.4514 | $\delta = 5.652905e-04$<br>itn = 1000<br>time sec = 108.8962 | $\delta = 6.498171e-04$<br>itn = 1000<br>time sec = 124.7879 | – |
| M = 10000 | $\delta = 2.472809e-03$<br>itn = 166<br>time sec = 4.6442 | $\delta = 9.357439e-04$<br>itn = 353<br>time sec = 18.0198 | $\delta = 1.605869e-03$<br>itn = 633<br>time sec = 51.0818 | $\delta = 2.498557e-03$<br>itn = 668<br>time sec = 59.8452 | $\delta = 2.851228e-03$<br>itn = 1000<br>time sec = 288.304 |

**Conclusions**. Results of the presented study show that penalty method (3) is fast since it uses the supplied program for calculation of subgradients but it requires careful selection of the penalty parameter M. Method (4) is slow since in the present study it uses finite differences for calculation of generalized gradients but is rather stable with respect the choice of the penalty parameter $M$ . Further studies will be directed on investigation of differential properties of the projection mapping $\pi_X(x)$ and the function $f(\pi_X(x))$ to reduce time of calculating subgradients of $f(\pi_X(x))$, also for account of parallel calculations. Some preliminary results in this direction are available in [15].

**Acknowledgment**. The research was supported by grant 2020.02/0121 of the National Research Fund of Ukraine.

**Vladimir Norkin,**
D.Sc. (fiz. and. math.), Leading researcher, V.M. Glushkov Institute of Cybernetics of the NAS of Ukraine, Kyiv,
https://orcid.org/0000-0003-3255-0405

**Anton Kozyriev,**
postgraduate student, National Technical University of Ukraine "Ihor Sikorsky Kyiv Polytechnic Institute", Kyiv.

**Vladimir Norkin [1, 2], Anton Kozyriev [2 *]**

## On Shor's *r*-Algorithm for Problems with Constraints

*[1] V.M. Glushkov Institute of Cybernetics of the NAS of Ukraine, Kyiv,*
*[2] National Technical University of Ukraine "Ihor Sikorsky Kyiv Polytechnic Institute", Kyiv*
*[*] Correspondence: a.kozyriev@kpi.ua*

**Introduction.** Nonsmooth optimization problems arise in a wide range of applications, including engineering, finance, and deep learning, where activation functions often have discontinuous derivatives, such as ReLU. Conventional optimization algorithms developed primarily for smooth problems face difficulties when applied to nonsmooth contexts due to discontinuities and other associated irregularities. Possible approaches to overcome these problems include smoothing of functions and applying non-smooth optimization techniques. In particular, Shor's r-algorithm (Shor, Zhurbenko (1971), Shor (1979)) with space stretching in the direction of the difference of two subsequent subgradients is a competitive non-smooth optimization method (Bagirov et al. (2014)). However, the original r-algorithm is designed to minimize unconstrained convex functions.

**The goal of the work.** The standard technique for applying this algorithm to problems with constraints consists in the use of exact non-smooth penalty functions (Eremin (1967), Zangwill (1967)). At the same time, it is necessary to correctly choose (quite large) the penalty coefficient of the penalty functions. Norkin (2020, 2022), Galvan et al. (2021) propose the so-called projective exact penalty functions method, which theoretically does not require choice of the penalty coefficient. The purpose of the present work is to study an applicability of the new exact projective non-smooth penalty functions method for solving conditional problems of non-smooth optimization by Shor's r-algorithm.

**The results.** In this paper, the original optimization problem with convex constraints is first transformed into an unconstrained problem by the projective penalty function method, and then the r-algorithm is used to solve the transformed problem. The results of testing this approach on problems with linear constraints using a program implemented in Matlab are presented. The results of the present study show that the standard method of non-smooth penalties combined with Shor's r-algorithm is fast, due to the use of the provided program to calculate the subgradients, but it requires the correct selection of the penalty parameter. The projective penalty method is slow because in this study it uses finite differences to calculate the gradients, but it is quite stable with respect to the choice of the penalty parameter. Further research will be aimed at investigating the differential properties of the projection mapping and reducing the time of computing subgradients for account of parallel calculations.

**В.І. Норкін [1,2], А.Ю. Козирєв [2*]**

# Про *r*-алгоритм Шора для задач з обмеженнями

[1] *Інститут кібернетики імені В.М. Глушкова НАН України, Київ*
[2] *Національний технічний університет України «Київський політехнічний інститут ім. Ігоря Сікорського»*
[*] *Листування:* *a.kozyriev@kpi.ua*

**Вступ.** Проблеми негладкої оптимізації виникають у широкому діапазоні застосувань, включаючи техніку, фінанси та глибоке навчання, де функції активації часто мають розривні похідні, наприклад, ReLU. Звичайні алгоритми оптимізації, розроблені переважно для гладких проблем, стикаються з труднощами при застосуванні до негладких контекстів через розриви та інші пов'язані нерегулярності. Можливі підходи для подолання цих проблем включають методи згладжування функцій та застосування методів негладкої оптимізації. Зокрема, r-алгоритм Шора (Шор, Журбенко (1971), Шор (1979)) з розтягуванням простору у напрямку різниці двох послідовних субградієнтів це конкурентно-здатний метод негладкої оптимізації (Bagirov et al. (2014)). Проте оригінальний r-алгоритм призначений для мінімізації опуклих яружних функцій без обмежень.

**Мета роботи.** Стандартний прийом для вирішення цим алгоритмом завдань з обмеженнями полягає у застосуванні точних негладких штрафних функцій (Єрьомін (1967), Zangwill (1967)). При цьому необхідно правильно вибрати (досить великий) коефіцієнт штрафу методом штрафних функцій. У роботах (Норкін (2020, 2022), Galvan та ін. (2021) запропоновано так званий проєктивний метод точних штрафних функцій, який теоретично не потребує точного визначення штрафного коефіцієнта. Мета роботи полягає у дослідженні можливостей нового методу точних проєктивних негладких штрафних функцій для розв'язання умовних задач негладкої оптимізації r-алгоритмом Шора.

**Результати.** У цій статті негладка оптимізаційна задача з опуклими обмеженнями спочатку перетворюється на задачу без обмежень методом проєктивних штрафних функцій, а потім для вирішення перетвореної задачі застосовується r-алгоритм. Наводяться результати тестування такого підходу на задачах з лінійними обмеженнями за допомогою програми, реалізованої у Matlab. Результати представленого дослідження показують, що стандартний метод негладких штрафів у поєднанні з r-алгоритмом Шора – швидкий, оскільки він використовує надану програму для розрахунку субградієнтів, але вимагає ретельного вибору параметра штрафу. Метод проєктивних штрафів – повільний, оскільки в цьому дослідженні він використовує кінцеві різниці для розрахунків узагальнених градієнтів, але досить стабільний щодо вибору параметра штрафу. Подальші дослідження будуть спрямовані на дослідження диференціальних властивостей проєкційного відображення та скорочення часу обчислення субградієнтів за рахунок паралельних обчислень.